\documentclass[final,12pt,nonatbib]{elsarticle}
{
\makeatletter
\let\c@author\relax
\makeatother
}

\usepackage{xcolor,soul,cancel}
\usepackage[color]{showkeys}

\usepackage{mathtools}
\definecolor{refkey}{rgb}{0,1,1}
\definecolor{labelkey}{rgb}{1,0,0}
\usepackage{amssymb}
\usepackage{amsmath}
\usepackage{amsthm}
\usepackage{graphicx}
\usepackage{breqn}
\usepackage{float}
\usepackage{xcolor}
\usepackage{cases}
\usepackage{color}
\usepackage{hyperref}
\usepackage{tabularx}
\usepackage{soul}
\usepackage{cancel}
\usepackage[normalem]{ulem}
\hypersetup{
    colorlinks=true, 
    linkcolor=blue, 
    urlcolor=red, 
    linktoc=all 
}
\journal{arXiv}

\newtheorem{thm}{Theorem}

\newtheorem{cor}{Corollary}
\newtheorem{prop}[thm]{Proposition}

\numberwithin{equation}{section}
\newcommand{\floor}[1]{\lfloor #1 \rfloor}
\newcommand{\eq} [1] {\begin{equation}\label{#1}\quad}
\newcommand{\en} {\end{equation}}
\newcommand{\scal}[1]{\langle#1\rangle}
\newcommand{\norm}[1]{\left\Vert#1\right\Vert}
\newcommand{\abs}[1]{\left\vert#1\right\vert}

\newcommand{\N}{\mathbb N}
\newcommand{\C}{\mathbb C}

\newcommand{\R}{\mathbb R}

\newcommand{\im}{\operatorname{Im}}

\newcommand{\re}{\operatorname{Re}}

\newcommand{\Rec}{\operatorname{\bf RM}}

\begin{document}

\begin{frontmatter}

\title{Unified approach to reciprocal matrices with Kippenhahn curves containing elliptical components }


\author[ucb]{Muyan Jiang}
\ead{muyan_jiang@berkeley.com, mj2259@nyu.edu, jimmymuyan@163.com}
\author[nyuad]{Ilya M. Spitkovsky}
\address[ucb]{Department of Industrial Engineering \& Operations Research, University of California Berkeley, 4141 Etcheverry Hall, Berkeley, CA 94720, USA}
\address[nyuad]{Mathematics Program, Division of Science and Mathematics, New York  University Abu Dhabi (NYUAD), Saadiyat Island,
P.O. Box 129188 Abu Dhabi, United Arab Emirates}
\ead{imspitkovsky@gmail.com, ims2@nyu.edu, ilya@math.wm.edu}
\tnotetext[support]{The second author [IMS]  was supported in part by Faculty Research funding from the Division of Science, New York University Abu Dhabi.}

\begin{abstract}
 {\em Reciprocal matrices}  are tridiagonal matrices $(a_{ij})_{i,j=1}^n$ with constant main diagonal and such that $a_{i,i+1}a_{i+1,i}=1$ for $i=1,\ldots,n-1$. For these matrices, criteria are established under which their Kippenhahn curves contain elliptical components or even consist completely of such. These criteria are in terms of system of homogeneous polynomial equations in variables $(\abs{a_{j,j+1}}-\abs{a_{j+1,j}})^2$, and established via a unified approach across arbitrary dimensions. The results are illustrated, and specific numerical examples provided, for $n=7$ thus generalizing earlier work in the lower dimensional setting.
\end{abstract}


\end{frontmatter}

\section{Introduction} 
As usual, let $\C^{n\times n}$ stand for the algebra of all $n$-by-$n$ matrices with complex entries, and for any $A\in \C^{n\times n}$ denote 
\[ \re A =\frac{A+A^*}{2}, \quad \im A=\frac{A-A^*}{2i}. \]
The characteristic polynomial $P_A(\lambda,\theta)=\det(\re(e^{i\theta}A)-\lambda I)$ of $\re(e^{i\theta}A)$ is sometimes called the {\em numerical range generating} polynomial, or the {\em Kippenhahn polynomial} of $A$. The reason behind this name is that the envelope $C(A)$ of the family of lines 
\[ \{e^{-i\theta}(\lambda_j(\theta)+i\R)\colon \theta\in (-\pi,\pi],\ j=1,\ldots,n\}, \]
where $\lambda_j(\theta)$ are the roots of $P_A(\lambda,\theta)$, is a curve the convex hull of which coincides with the {\em numerical range} 
\[ W(A):= \{ \scal{Ax,x}\colon x\in\C^n, \norm{x}=1\} \] 
of $A$ --- the fact established by Kippenhahn in \cite{Ki} (see also the English translation \cite{Ki08}). Respectively, $C(A)$ is called the {\em numerical range generating curve}, or the {\em Kippenhahn curve}, of $A$. 
Moreover, as can be seen from \cite{Woe08}, $C(A)$ defines completely the so called {\em rank-$k$ numerical ranges}  
\[ \Lambda_k(A)=\{ \lambda\in\C\colon PAP=\lambda P \text{ for some orthogonal rank } k \text{ projection } P\} \]
of $A$ for all values of $k$, not only for $k=1$ when $\Lambda_1(A)=W(A)$.

Note also that the foci of $C(A)$ are exactly the eigenvalues of $A$. 

In this paper, we continue studying Kippenhahn curves of {\em reciprocal} matrices -- the class introduced (to the best of our knowledge) in \cite{BPSV} and treated further in \cite{BPS, JiangS}. Remind therefore that a matrix $A=(a_{ij})_{i,j=1}^n$ is reciprocal if it is {\em tridiagonal} (i.e., $a_{ij}=0$ if $\abs{i-j}>1$), has constant main diagonal and, in addition, $a_{i,i+1}a_{i+1,i}=1$ for $i=1,\ldots,n-1$. 

All the objects we are interested in ($P_A(\lambda,\theta), C(A), W(A)$) behave in a natural and predictable way under translations of any $A\in\C^{n\times n}$. For reciprocal matrices we therefore may (and will) without loss of generality suppose that they have the zero main diagonal. The class of $n$-by-$n$ reciprocal matrices satisfying this additional assumption in what follows will be denoted $\Rec_n$. 

By Corollary 2 in \cite{BPSV}, the Kippenhahn curve \(C(A)\) of a matrix $A\in \Rec_n$ is symmetric about both coordinate axes. Moreover, all these matrices have the same spectrum
\eq{spectrum}
\sigma(A) = \left\{2\cos\frac{k \pi}{n+1}\colon k=1,\cdots,n\right\}.
\en
So, if $C(A)$ contains an ellipse $E$ then the latter has real foci, and it is either centered at the origin or its reflection $-E$ about the $y$-axis is also contained in $C(A)$. In the latter case, we will be saying that $C(A)$ contains a pair of shifted ellipses $E\cup(-E)$. 

We aim to develop a generic procedure to deduce criteria for the following three cases: when a Kippenhahn curve (i) contains an ellipse centered at the origin, (ii) consists of concentric ellipses all centered at the origin, and (iii) contains a pair of shifted ellipses. In each case, the criterion is a system of algebraic equations, most conveniently written in terms of \eq{xij} \xi_j=(\abs{a_{j,j+1}}-\abs{a_{j+1,j}})^2/4, \quad j=1,\ldots, n-1. \en While we standardize the criteria derivation process, we make quantitative observations of the respective system.

These criteria are obtained in Section~\ref{s:theory}, and implemented for $n=7$ in Sections~\ref{s:n=7}--\ref{s:shif7ex}. More specifically, Section~\ref{s:n=7} provides the explicit conditions for $A\in\Rec_7$ to have $C(A)$ containing at least one (Subsection~\ref{s:onel7}) or all three (Subsection~\ref{s:conel7}) ellipses centered at the origin. In Section~\ref{s:shif7}, it is established that $C(A)$ may contain a pair of shifted ellipses $E\cup(-E)$ only if the foci of $E$ (and thus $-E$ as well) have opposite signs, confirming the conjecture of \cite{JiangS} in the case $n=7$. Consequently, there are exactly three possible configurations of $C(A)$ containing a pair of shifted ellipses. They are described completely in Section~\ref{s:shif7ex}. An auxiliary Section~\ref{s:prel} contains some technical formulas for the Kippenhahn polynomial of reciprocal matrices.
\section{Preliminaries} \label{s:prel} 
The following technical observation will prove itself useful: the Kippenhahn curve $C(A)$ of any (not necessarily reciprocal) matrix $A$ contains an ellipse $E$ with real foci if and only if its Kippenhahn polynomial is divisible by 
\eq{polE} (\lambda-p\cos\theta)^2-(X^2\rho+c^2). \en
Here $p$ is the center of $E$ while $c$ and $X$ stand for the length of the minor half-axis and half the distance between the foci, respectively (see, e.g., Proposition~1 in \cite{JiangS}).  

We turn therefore to the specifics of the Kippenhahn polynomial of reciprocal matrices.

It was observed in \cite{BPSV} that the Kippenhahn polynomial of $A\in\Rec_n$ has the form 
\eq{Pn} P_n(\zeta)=\zeta^m+\sum_{j=0}^{m-1}p_j\zeta^j, \en
pre-multiplied by $-\lambda$ if $n$ is odd. Here $\zeta=\lambda^2, m=\floor{n/2}$, $p_j$ depend only on $A_k:=(\abs{a_{k,k+1}}^2+\abs{a_{k+1,k}}^2)/2$, $k=1,\ldots,n-1$, and $\tau:=\cos(2\theta)$. The respective explicit formulas were obtained and used there for $n=4,5,6$. 

As was figured out in \cite{JiangS}, for these values of $n$ the coefficients $p_j$ in \eqref{Pn} take a simpler form when expressed in terms of $\xi_k$ given by\eqref{xij} in place of $A_k$, and $\rho=\cos^2\theta$ in place of $\tau$. Here we take the matter a step further, obtaining explicit formulas for $p_j$ for arbitrary values of $n$. 

In order to achieve this, we need an explicit formula for determinants of tridiagonal matrices with a constant main diagonal. It is convenient to introduce the notation $i\prec j$ if $i,j\in\N$ and $j-i>1$.

\begin{prop}\label{th:exptri}Let $X=(x_{ij})$ be a tridiagonal $n$-by-$n$ matrix. If $x_{11}=\ldots=x_{nn}=:a$, then 
\[ \det X = \begin{cases} D_n & \text{ if } n \text{ is even,} \\ aD_n & \text{ if } n \text{ is odd,} \end{cases} \]
with \eq{Dn} D_n=a^{2m}+\sum_{j=0}^{m-1}d_ja^{2j}.\en 
Here $m=\lfloor{n/2}\rfloor$, 
\eq{dj} d_{m-j}= (-1)^j\sum_{1\leq i_1\prec i_2\prec\cdots\prec i_j<n}\prod_{k=1}^j \eta_{i_k}, \quad j=1,\ldots,m,\en 
while $\eta_i:=x_{i,i+1}x_{i+1,i}$.  In particular, each coefficient $d_j$ is a homogeneous polynomial of degree $m-j$ in  $\eta_1,\ldots,\eta_{n-1}$, making $P_n$ homogeneous of degree $m$ in  $\eta_1\ldots,\eta_{n-1}$ and $a^2$.
\end{prop} 
\begin{proof}Recall that for any tridiagonal matrix $X$, whether or not its main diagonal is constant, denoting by $\Delta_j$ the determinant of its left upper $j$-by-$j$: 
\eq{recur} \Delta_n=x_{nn}\Delta_{n-1}-x_{n-1,n}x_{n,n-1}\Delta_{n-2}. \en 
Using this well-known recursive relation and taking into consideration the definition of $D_n$:  \eq{Dnrec} D_n=\begin{cases} D_{n-1}-\eta_{n-1}D_{n-2},  \\ a^2 D_{n-1}-\eta_{n-1}D_{n-2}, \end{cases}  \text{ if } n \text { is } \begin{cases} \text{odd}, \\ \text{ even}. \end{cases} \en
Relabel temporarily $d_j$ from \eqref{Dn} by $d_j^{(n)}$, to emphasize its dependence on $n$. In this notation, \eqref{Dnrec} implies 
\eq{dnrec} d^{(n)}_{\floor{n/2}-j}=d^{(n-1)}_{\floor{(n-1)/2}-j}-\eta_{n-1}\,
d^{(n-2)}_{\floor{(n-2)/2}-(j-1)}.\en
Since $D_1=1$ and $D_2=a^2-\eta_1$, formulas \eqref{dj} hold in a trivial way for $n=1,2$. The recursion \eqref{dnrec}, combined with the mathematical induction principle, show the validity of \eqref{dj} for all $n$. 
Indeed, under the induction hypotheses the right-hand side of \eqref{dnrec} is 
\begin{align*} & (-1)^j\sum_{1\leq i_1\prec i_2\prec\cdots\prec i_j<n-1}\prod_{k=1}^j \eta_{i_k} -(-1)^{j-1}\sum_{1\leq i_1\prec i_2\prec\cdots\prec i_{j-1}<n-2}\eta_{n-1}\prod_{k=1}^{j-1} \eta_{i_k} \\ & =
 (-1)^j\left(\sum_{1\leq i_1\prec i_2\prec\cdots\prec i_j<n-1}\prod_{k=1}^j \eta_{i_k}+ \sum_{1\leq i_1\prec i_2\prec\cdots\prec i_{j}=n-1}\prod_{k=1}^{j}\eta_{i_k}\right) \\
& =  (-1)^j\sum_{1\leq i_1\prec i_2\prec\cdots\prec i_j<n}\prod_{k=1}^j \eta_{i_k}.
 \end{align*} \end{proof} 

For any $A\in\Rec_n$, the matrix $\re(e^{i\theta}A)-\lambda I$ is tridiagonal along with $A$, all its diagonal entries equal $-\lambda$, and the product of the off-diagonal entries in the $(i,i+1)$ and $(i+1,i)$ positions is 
\[ \frac{1}{4}\left(e^{i\theta}a_{i,i+1}+e^{-i\theta}\overline{a_{i,i+1}}\right)\left(e^{-i\theta}\overline{a_{i+1,i}}+e^{i\theta}a_{i+1,i}\right)=\xi_i+\rho. \]
Invoking Proposition~\ref{th:exptri}, we immediately arrive at 
\begin{prop} \label{th:expf} For any $A\in\Rec_n$ the coefficients of the polynomial \eqref{Pn} are given by \eq{pj} p_{m-j}= (-1)^j\sum_{1\leq i_1\prec i_2\prec\cdots\prec i_j<n}\prod_{k=1}^j (\xi_{i_k}+\rho), \quad j=1,\ldots,m.\en \end{prop} 
On the other hand, when applying Proposition~\ref{th:exptri} to $A-\lambda I$ we have $\eta_1=\ldots=\eta_{n-1}=1$, and so the right-hand side of \eqref{dj} is $(-1)^j$ times the number of non-consecutive selection 
of $j$ objects out of $n-1$. Since the latter is nothing but $n-j\choose j$, the characteristic polynomial of $A$ is \eq{charpol} Q_n(\lambda)= \sum_{j=0}^m (-1)^j{n-j\choose j}\lambda^{2(m-j)},\en 
premulitplied by $-\lambda$ in case of odd $n$. 

\section{Elliptical components of the Kippenhahn curve} \label{s:theory}
\subsection{An ellipse centered at the origin} \label{s:onel}

For $E$ centered at the origin \eqref{polE} takes the form $\zeta-(C+\rho X^2)$, where, along with the abbreviation $\zeta=\lambda^2$ we have also relabeled $C:=c^2$. 

By Bezout's theorem, the divisibility of $P_n$ by $\zeta-(C+\rho X^2)$ is equivalent to 
\eq{eq0} P_n(C+\rho  X^2) = 0.\en
According to \eqref{Pn} and \eqref{pj}, the left-hand side of \eqref{eq0} can be rewritten as a polynomial in $\rho$ of degree $m$ with the coefficients $f_j(\xi,C,X)$ of $\rho^j$ being polynomials  in $X,C,\xi_1,\ldots,\xi_{n-1}$. For \eqref{eq0} to hold identically in $\rho$, it is necessary and sufficient that 
\eq{fxic} f_j(\xi,C,X)=0, \quad j=0,\ldots,m. \en 

Note that the polynomials $f_j$ are homogeneous of degree $m-j$ in $C$ and $\xi_1,\ldots,\xi_{n-1}$, if $X$ is treated as a parameter. In particular, $f_m$ depends on $X$ only. 

Direct computations show that $f_m(\xi,C,X)=Q_n(X)$, and is therefore equal to zero if and only if $X$ is as given by \eqref{spectrum}. From its positivity follows that only $k=1,\ldots,m$ are admissible. 

Next, let us use the fact that $f_{m-1}$ is linear as a function of $C,\xi_1,\ldots,\xi_{n-1}$. Moreover, the coefficient of $C$ is 
\[ \sum_{j=0}^{m-1}(-1)^j(m-j){n-j\choose j}X^{2(m-j-1)},\]
which (up to a multiple $2X$) is $Q_n'(X)$. Since the eigenvalues of $A$ are simple, so are the roots of $Q_n$. Consequently, for any fixed value of $X=2\cos\frac{k\pi}{n+1}$ the $m$th equation in \eqref{fxic} is satisfied while the ($m-1$)st can be used to express $C$ uniquely as a linear combination of $\xi_1,\ldots,\xi_{n-1}$. 

Plugging the resulting formula for $C$ in \eqref{fxic} with $j=0,1,\ldots,m-2$, we arrive at the system of $m-1$ polynomial equations homogeneous in $\xi_1,\ldots,\xi_{n-1}$ and having degrees ranging from 2 to $m$. 

Therefore, the following result holds. 
\begin{thm}\label{th:crionel}Let $A=\left(a_{ij}\right)\in \Rec_n$. Then for $C(A)$ to contain an ellipse $E$ centered at the origin the following conditions are necessary and sufficient: 

{\em (i)} the foci of $E$ are located at $\pm X$ for $X=2\cos\frac{k\pi}{n+1}$ and some $k=1,\ldots,m:=\lfloor\frac{n}{2}\rfloor$; 

{\em (ii)} Defined by \eqref{xij} variables $\xi_j$ satisfy the system of $m-1$ homogeneous equations obtained as described above. \end{thm} 

\subsection{Concentric Ellipses}\label{s:conel} 

Based on Theorem~\ref{th:crionel}, it is straightforward to state the criteria for $C(A)$ to contain several ellipses centered at the origin with the prescribed lengths of their major axes. For $k$ ellipses, the respective criterion will consist of $k(m-1)$ homogeneous polynomial equations in $\xi_j$. For $k=m$, however, just $(m-1)^2$ equations would suffice. Indeed, if $C(A)$ already contains $m-1$ elliptical components, the remaining one is forced to be an ellipse as well. 

Here is an alternative approach to the case of $m$ concentric ellipses,  which allows for a further reduction in the number of equations.

Namely, for this case to materialize it is necessary and sufficient that the polynomial $P_n$ as in \eqref{Pn} is divisible by $\zeta-(C_k+\rho X_k^2)$ for some $C_k>0$ and $X_k=2\cos\frac{k\pi}{n+1}$ ($k=1,\ldots,m$), and therefore by their product
\eq{prod}
\prod _{k=1}^m \left(\zeta -\left(C_k+\rho X_k^2\right)\right).
\en
Since the highest term in both \eqref{Pn} and \eqref{prod} is $\zeta^m$, the two polynomials coincide. With the use of \eqref{pj}, equating the respective coefficients of \eqref{Pn} and \eqref{prod} yields  
\eq{coeff}  \sum_{1\leq i_1<\ldots<i_j\leq m}\prod_{k=1}^j\left(C_{i_k}+\rho X_{i_k}^2\right) = \sum_{1\leq i_1\prec i_2\prec\cdots\prec i_j<n}\prod_{k=1}^j (\xi_{i_k}+\rho), \en
$j=1,\ldots,m$. 

In its turn, both sides of \eqref{coeff} are polynomials of degree $j$ in $\rho$. The leading coefficients of these polynomials coincide automatically, since 
\[ \sum_{1\leq i_1<\ldots<i_j\leq m}\prod_{k=1}^j X_{i_k}^2 = {n-j \choose j}. \] 

Equating the coefficients of $\rho^i$ in the $j$th condition in \eqref{coeff} for the remaining values of $i=0,\ldots,j-1$, 
we arrive at $j$ equations, respectively, homogeneous of degree $j-i$ in $C_k,\xi_k$. The total number of these equations is therefore $m(m+1)/2$. 

Let us tackle first $m$ equations corresponding to $i=j-1$. They can be rewritten as a system of linear equations in $C_1,\ldots, C_m$ with the right-hand sides being linear functions of $\xi_1,\ldots,\xi_{n-1}$. The matrix $Z$ of this system is given by
\eq{matrix} z_{pq}=\sum_{i_1<\ldots < i_{q-1}, i_k\neq p}\prod_{k=1}^{q-1} X_{i_k}^2, \quad p,q=1,\ldots m. \en
In particular, the first row of $Z$ consists of all ones.

Since $\det Z=\prod_{i,j=1,\ldots,m;i>j}(X_i^2-X_j^2)\neq 0$, we can represent $C_1,\ldots,C_m$ uniquely as linear functions of $\xi_1,\ldots,\xi_{n-1}$. 

Substituting these expressions for $C_k$ into the remaining $m(m-1)/2$ equations, we arrive at the system containing only $\xi_1,\ldots,\xi_{n-1}$ as the unknowns. More specifically, it is comprised of $k$ equations of degree $m-k+1$, with $k$ running from 1 through $m-1$. 

To summarize: 
\begin{thm}\label{th:criallel}Let $A=\left(a_{ij}\right)\in \Rec_n$. Then for $C(A)$ to consist of concentric ellipses it is necessary and sufficient that $\xi_j:=(\abs{a_{j,j+1}}-\abs{a_{j+1,j}})^2/4$ ($j=1,\ldots,n-1$) satisfy the system of $\lfloor n\rfloor \left(\lfloor n\rfloor-1\right)/2$ homogeneous polynomial equations obtained as described above. \end{thm} 

More specifically, under the conditions of Theorem~~\ref{th:criallel}  $C(A)=\cup_{j=1}^m E_j$, where $E_j$ is an ellipse with \sout{the} foci $\pm 2\cos\frac{j\pi}{n+1}$ and \sout{the} minor axis of length $2\sqrt{C_j}$. For any admissible ($n-1$)-tuple ${\xi_1,\ldots,\xi_{n-1}}$ the respective values of $C_j$ can be found from the system of linear equations with the $Z$ as the matrix of coefficients. Alternatively, they are nothing but the roots (in decreasing order) of \eqref{Pn} with $p_j$ defined by \eqref{pj} for $\rho=0$.

\subsection{Shifted Ellipses} \label{s:shifel} 

Recall that the pair of shifted ellipses, if contained in the Kippenhahn curve $C(A)$ of a matrix $A\in\Rec_n$, has to be of the form $E\cup (-E)$. Without loss of generality, let $E$ be the ellipse centered at $p>0$. Then the ellipse $-E$ is centered at $-p$ while the length $c$ of the minor half-axis and half the distance $X$ between the foci are the same for $E$ and $-E$. Invoking \eqref{polE}, we conclude that \eq{pash}C(A)\supset E\cup (-E) \text{ with the ellipse } E \text{ different from } -E  \en  if and only if the polynomial \eqref{Pn} is divisible by
\begin{multline}\label{shell}  \left((\lambda-p\cos\theta)^2-(X^2\rho+c^2)\right)\left((\lambda+p\cos\theta)^2-(X^2\rho+c^2)\right) \\ = \zeta ^2-2\zeta  \left(C+\rho  \left(X^2+p^2\right)\right)+\left(C+\rho  \left(X^2-p^2\right)\right)^2. \end{multline}
Here we are again using the abbreviation $C=c^2$, along with $\lambda^2=\zeta$. 

Equivalently, 
\[ P_n(\zeta_+)=P_n(\zeta_-)=0 \text{ identically in } \rho, \] 
where 
\[ \zeta_{\pm} = C + (p^2 + X^2) \rho \pm 2p\sqrt{\rho(C+X^2\rho)} \]
are the roots of \eqref{shell} considered as the quadratic polynomial in $\zeta$.  

Observe that \eq{roe} R_e:=P_n(\zeta_+)+P_n(\zeta_-) \text{ and } R_o:=\frac{P_n(\zeta_+)-P_n(\zeta_-)}{\sqrt{\rho(C+X^2\rho)}} \en  are polynomials in $\rho$ of degrees at most $m$ and $m-1$, respectively. 

Note, in addition, that asymptotically 
\[ \zeta_\pm \sim (X\pm p)^2\rho=f^2_\pm\rho, \]
where $f_\pm:=p\pm X$ are the foci of $E$. 

From \eqref{Pn} and \eqref{pj} therefore: 

\[ P_n(\zeta_\pm) \sim \rho^m \sum_{j=0}^m (-1)^j{n-j \choose j} f_\pm^{2(m-j)}=Q_n(f_\pm^2), \]
where $Q_n$ is given by \eqref{charpol}. Consequently, the leading coefficients of $R_e, R_o$ vanish if and only if $f_\pm$ are (non-zero) eigenvalues of $A$ --- the condition already imposed on $X,p$ earlier. 

There are $m(m-1)$ pairs ${X,p}$ satisfying this condition. For each such pair, considered as parameters, the coefficients of $\rho^j$ in $R_e$ and $R_o$ are homogeneous polynomials in $C,\xi_1,\ldots,\xi_{n-1}$ of degree $m-j$ and $m-j-1$, respectively.

Equating these $2m-1$ coefficients with zero, use the one corresponding to $\rho^{m-1}$ in $R_e$ (or $\rho^{m-2}$ in $R_o$) to express $C$ as a linear function of $\xi_1,\ldots,\xi_{n-1}$. Substituting this expression into the remaining $2m-2$ equations, we arrive at the criterion for \eqref{pash} to hold.

\begin{thm}\label{th:shifel}Let $A=\left(a_{ij}\right)\in \Rec_n$. Then for \eqref{pash} to hold it is necessary and sufficient that $\xi_j:=(\abs{a_{j,j+1}}-\abs{a_{j+1,j}})^2/4$ ($j=1,\ldots,n-1$) satisfy the system of $2m-2$ homogeneous polynomial equations obtained as described above. \end{thm} 

\section{Origin-centered ellipses in $n=7$ case} \label{s:n=7}

To illustrate the approach developed in Section~\ref{s:theory}, consider matrices $A\in \Rec_7$. This is a natural step forward, given the constructive results for reciprocal matrices of sizes up to 6-by-6 obtained in \cite{BPSV,JiangS}.

As is the case for any odd $n$, $C(A)$ contains the origin. Disregarding this, since $m=\lfloor 7/2\rfloor =3$, {\em a priori} possibilities are exactly the same as for $n=6$: the Kippenhahn curve $C(A)$ has no elliptical components, contains exactly one ellipse (which is then centered at the origin), consists of three concentric ellipses (all centered at the origin), or of one ellipse centered at the origin and two shifted ellipses.

Applying \eqref{pj} to the case $n=7$, from \eqref{Pn} we obtain:

\begin{dmath}\label{P7}
P_7(\zeta)= \zeta ^3+\zeta ^2 \left(-\xi _1-\xi _2-\xi _3-\xi _4-\xi _5-\xi _6-6 \rho \right)+\zeta  \left(\left(4 \xi _1+3 \xi _2+3 \xi _3+3 \xi _4+3 \xi _5+4 \xi _6\right) \rho +\xi _1 \xi _3+\xi _1 \xi _4+\xi _2 \xi _4+\xi _1 \xi _5+\xi _2 \xi _5+\xi _3 \xi _5+\xi _1 \xi _6+\xi _2 \xi _6+\xi _3 \xi _6+\xi _4 \xi _6+10 \rho ^2\right)-4 \rho ^3+\left(-3 \xi _1-\xi _2-2 \xi _3-2 \xi _4-\xi _5-3 \xi _6\right) \rho ^2+\left(-2 \xi _1 \xi _3-\xi _5 \xi _3-\xi _6 \xi _3-\xi _1 \xi _4-\xi _2 \xi _4-\xi _1 \xi _5-2 \xi _1 \xi _6-\xi _2 \xi _6-2 \xi _4 \xi _6\right) \rho -\xi _1 \xi _3 \xi _5-\xi _1 \xi _3 \xi _6-\xi _1 \xi _4 \xi _6-\xi _2 \xi _4 \xi _6. \end{dmath}

\subsection{One ellipse} \label{s:onel7}
Equations \eqref{fxic} for $j=0,1,2$ take the form 
\begin{dmath}\label{fxic73}
    C^3-C^2 \left(\xi _1+\xi _2+\xi _3+\xi _4+\xi _5+\xi _6\right)+C \left(\xi _1 \xi _3+\xi _5 \xi _3+\xi _6 \xi _3+\xi _1 \xi _4+\xi _2 \xi _4+\xi _1 \xi _5+\xi _2 \xi _5+\xi _1 \xi _6+\xi _2 \xi _6+\xi _4 \xi _6\right)-\left(\xi _1 \xi _3 \xi _5+\xi _1 \xi _3 \xi _6+\xi _1 \xi _4 \xi _6+\xi _2 \xi _4 \xi _6\right)=0,
\end{dmath}
\begin{dmath}\label{fxic72} 
    3 C^2 \left(X^2-2\right)-C \left(2 \left(X^2-2\right)\xi _1+ \left(2 X^2-3\right)\left(\xi _2+\xi _3+\xi _4+\xi _5\right)+2\left(X^2-2\right)\xi _6\right)+X^2\xi _2 \xi _5 +\left(X^2-2\right)\left(\xi _1 \xi _3+\xi _1 \xi _6+\xi _4 \xi _6\right) +\left(X^2-1\right)\left(\xi _1 \xi _4+\xi _2 \xi _4+\xi _1 \xi _5+\xi _3 \xi _5+\xi _2 \xi _6+\xi _3 \xi _6\right)=0,
\end{dmath} and
\begin{dmath}
\label{fxic7}
    C \left(3 X^4-12 X^2+10\right)- \left(X^4-3 X^2+2\right)\left(\xi _3+\xi _4\right)-\left(X^4-3 X^2+1\right)\left(\xi _2+\xi _5\right) - \left(X^4-4 X^2+3\right)\left(\xi _1+\xi _6\right)=0,
\end{dmath}
respectively, with the parameter $X$ assuming one of the values \[ X_1=2\cos\frac{\pi}{8}=\sqrt{2+\sqrt{2}},   
X_2=2\cos\frac{\pi}{4}=\sqrt{2},  X_3=2\cos\frac{3\pi}{8}=\sqrt{2-\sqrt{2}}.  \] 
In agreement with the observation made in Section~\ref{s:onel}, the coefficient of $C$ in \eqref{fxic7} is nonzero for all the choices of $X$ and can thus be used to rewrite \eqref{fxic73}, \eqref{fxic72} as the system of two homogeneous equations for $\xi_1,\ldots,\xi_6$, one quadratic and one cubic. This is the realization of Theorem~\ref{th:crionel} for the case of 7-by-7 matrices.

Here are the numerical specifics for the case of $X=X_2$; the other two numerical values can be treated similarly.

As it happens, \eqref{fxic7} boils down to $C = \frac{1}{2} \left(\xi _1+\xi _2+\xi _5+\xi _6\right)$. Plugging this into \eqref{fxic73} and \eqref{fxic72}, we get
\[
\begin{cases}
\begin{aligned}
-\frac{1}{8} \bigl(&\xi _1+\xi _2-\xi _5-\xi _6\bigr) \bigl(\xi _1^2+2 \left(\xi _2-\xi _3-\xi _4\right) \xi _1-\xi _6^2 \\
&+\left(\xi _2+2 \xi _3-2 \xi _4-\xi _5\right) \left(\xi _2+\xi _5\right)+2 \left(\xi _3+\xi _4-\xi _5\right) \xi _6\bigr) & = 0,
\end{aligned}   \\
\frac{1}{2} \left(-\xi _2-\xi _3+\xi _4+\xi _5\right) \left(\xi _1+\xi _2-\xi _5-\xi _6\right)  = 0. & 
\end{cases}
\]
Consequently, the following result holds.
\begin{cor}\label{th:midel7}Let $A\in\Rec_7$. Then $C(A)$ contains an ellipse with the foci $\pm\sqrt{2}$
if and only if either 
\eq{crit71} \xi _1+\xi _2=\xi _5+\xi _6\en
or
\[ \begin{cases}
    \begin{aligned}
      \xi _1^2+2 \xi _2 \xi _1-2 \xi _3 \xi _1-2 \xi _4 \xi _1+\xi _2^2-\xi _5^2-\xi _6^2+2 \xi _2 \xi _3\\-2 \xi _2 \xi _4+2 \xi _3 \xi _5-2 \xi _4 \xi _5+2 \xi _3 \xi _6+2 \xi _4 \xi _6-2 \xi _5 \xi _6=0,
    \end{aligned}  \\
    \xi _2+\xi _3=\xi _4+\xi _5. 
\end{cases}
\] 
\end{cor} 

\begin{figure}[H]
    \centering
    \includegraphics[width=0.5\textwidth]{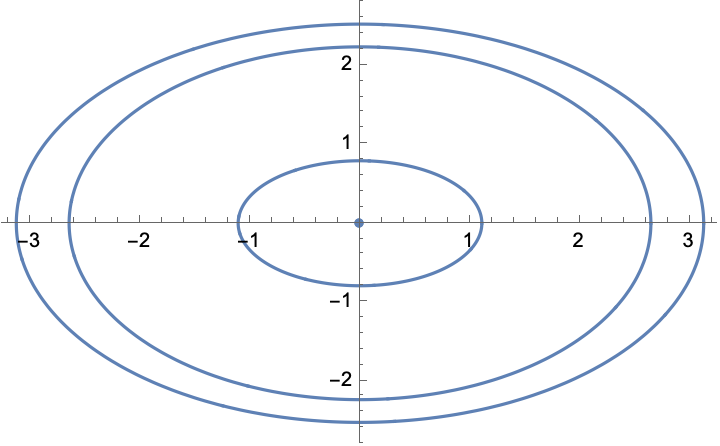}
    \caption{$\vec{\xi}$ = \{1, 4, 1, 1, 2, 3\}.}
\label{fig:1} 
\end{figure}

Let, for example,
\eq{ex71} \xi_1=\xi_3=\xi_4=1, \xi_2=4, \xi_5=2, \xi_6 = 3. \en 
Since \eqref{crit71} holds, the middle component of the Kippenhahn curve of the respective matrix $A$ is elliptical. See Fig.~\ref{fig:1} for the plot of $C(A)$ in this case. 

\subsection{Concentric ellipses} \label{s:conel7} 
We will follow the procedure described in Section~\ref{s:conel}. Equating the respective coefficients in the left- and right-hand sides of \eqref{coeff} considered as polynomials in $\rho$, we arrive at the following system of six equations in $C_j,\xi_j$:

\begin{align*}\label{sixfor7}
    &C_1+C_2+C_3=\xi _1+\xi _2+\xi _3+\xi _4+\xi _5+\xi _6,\\&
    \left(4-\sqrt{2}\right) C_1+4 C_2+\left(\sqrt{2}+4\right) C_3=4 \xi _1+3 \xi _2+3 \xi _3+3 \xi _4+3 \xi _5+4 \xi _6, \\&
    2 \left(2-\sqrt{2}\right) C_1+2 C_2+2 \left(\sqrt{2}+2\right) C_3=3 \xi _1+\xi _2+2 \xi _3+2 \xi _4+\xi _5+3 \xi _6, \\&
    C_1 C_2+C_3 C_2+C_1 C_3=\xi _1 \xi _3+\xi _5 \xi _3+\xi _6 \xi _3+\xi _1 \xi _4+\xi _2 \xi _4+\xi _1 \xi _5+\xi _2 \xi _5\\&+\xi _1 \xi _6+\xi _2 \xi _6+\xi _4 \xi _6, \\&
     \left(2-\sqrt{2}\right) C_1 C_2+\left(\sqrt{2}+2\right) C_3 C_2+2 C_1 C_3\\&=2 \xi _1 \xi _3+\xi _5 \xi _3+\xi _6 \xi _3+\xi _1 \xi _4+\xi _2 \xi _4+\xi _1 \xi _5+2 \xi _1 \xi _6+\xi _2 \xi _6+2 \xi _4 \xi _6, \\&
     C_1 C_2 C_3=\xi _1 \xi _3 \xi _5+\xi _1 \xi _3 \xi _6+\xi _1 \xi_4 \xi _6+\xi _2 \xi _4 \xi _6.\\
\end{align*}
The first three of these equations are linear; considering $C_j$ as the unknowns, the respective matrix $Z$ of coefficients is 
\begin{equation*}
    \left[
\begin{array}{ccc}
 1 & 1 & 1 \\
 4-\sqrt{2} & 4 & 4+\sqrt{2} \\
 4-2 \sqrt{2} & 2 & 4 +2 \sqrt{2} \\
\end{array}
\right].
\end{equation*}
This agrees with \eqref{matrix}, according to which for $n=7$ the matrix $Z$ should equal
\[
\left[
\begin{array}{ccc}
 1 & 1 & 1 \\
X_2 + X_3 & X_1 + X_3 & X_1 + X_2 \\
X_2 X_3 & X_1 X_3 & X_1 X_2 \\
\end{array}
\right]. 
\]

Solving these equations for $C_1,C_2,C_3$ in terms of $\xi_1,\ldots,\xi_6$: 
\begin{align*}
    &C_1 = \frac{1}{4} \left(\xi _1+\left(\sqrt{2}+1\right) \xi _2+\left(\sqrt{2}+2\right) \xi _3+\left(\sqrt{2}+2\right) \xi _4+\left(\sqrt{2}+1\right) \xi _5+\xi _6\right),\\
    &C_2 =\frac{1}{2} \left(\xi _1+\xi _2+\xi _5+\xi _6\right),\\
    &C_3 =  \frac{1}{4} \left(\xi _1+\left(1-\sqrt{2}\right) \xi _2+\left(2-\sqrt{2}\right) \xi _3+\left(2-\sqrt{2}\right) \xi _4+\left(1-\sqrt{2}\right)\xi _5+\xi _6\right),
\end{align*}
and plugging the result into the remaining three equations: 
\begin{dmath*}
\xi_1^3 + 3 \xi_2 \xi_1^2 + 4 \xi_3 \xi_1^2 + 4 \xi_4 \xi_1^2 + 3 \xi_5 \xi_1^2 + 3 \xi_6 \xi_1^2 + \xi_2^2 \xi_1 + 2 \xi_3^2 \xi_1 + 2 \xi_4^2 \xi_1 + \xi_5^2 \xi_1 + 3 \xi_6^2 \xi_1 + 4 \xi_2 \xi_3 \xi_1 + 4 \xi_2 \xi_4 \xi_1 + 4 \xi_3 \xi_4 \xi_1 + 2 \xi_2 \xi_5 \xi_1 - 28 \xi_3 \xi_5 \xi_1 + 4 \xi_4 \xi_5 \xi_1 + 6 \xi_2 \xi_6 \xi_1 - 24 \xi_3 \xi_6 \xi_1 - 24 \xi_4 \xi_6 \xi_1 + 6 \xi_5 \xi_6 \xi_1 - \xi_2^3 - \xi_5^3 + \xi_6^3 + 2 \xi_2 \xi_3^2 + 2 \xi_2 \xi_4^2 - 3 \xi_2 \xi_5^2 + 3 \xi_2 \xi_6^2 + 4 \xi_3 \xi_6^2 + 4 \xi_4 \xi_6^2 + 3 \xi_5 \xi_6^2 + 4 \xi_2 \xi_3 \xi_4 - 3 \xi_2^2 \xi_5 + 2 \xi_3^2 \xi_5 + 2 \xi_4^2 \xi_5 + 4 \xi_3 \xi_4 \xi_5 + \xi_2^2 \xi_6 + 2 \xi_3^2 \xi_6 + 2 \xi_4^2 \xi_6 + \xi_5^2 \xi_6 + 4 \xi_2 \xi_3 \xi_6 - 28 \xi_2 \xi_4 \xi_6 + 4 \xi_3 \xi_4 \xi_6 + 2 \xi_2 \xi_5 \xi_6 + 4 \xi_3 \xi_5 \xi_6 + 4 \xi_4 \xi_5 \xi_6 = 0,
\end{dmath*}
\begin{dmath*}
5 \xi_1^2 + 10 \xi_2 \xi_1 - 4 \xi_3 \xi_1 - 4 \xi_4 \xi_1 - 6 \xi_5 \xi_1 - 6 \xi_6 \xi_1 + 3 \xi_2^2 + 2 \xi_3^2 + 2 \xi_4^2 + 3 \xi_5^2 + 5 \xi_6^2 + 8 \xi_2 \xi_3 - 8 \xi_2 \xi_4 + 4 \xi_3 \xi_4 - 10 \xi_2 \xi_5 - 8 \xi_3 \xi_5 + 8 \xi_4 \xi_5 - 6 \xi_2 \xi_6 - 4 \xi_3 \xi_6 - 4 \xi_4 \xi_6 + 10 \xi_5 \xi_6 = 0,
\end{dmath*}
\begin{dmath}\label{con7}
5 \xi_1^2 + 6 \xi_2 \xi_1 - 8 \xi_3 \xi_1 - 2 \xi_5 \xi_1 - 6 \xi_6 \xi_1 - \xi_2^2 + 2 \xi_3^2 + 2 \xi_4^2 - \xi_5^2 + 5 \xi_6^2 + 4 \xi_2 \xi_3 - 4 \xi_2 \xi_4 + 4 \xi_3 \xi_4 - 2 \xi_2 \xi_5 - 4 \xi_3 \xi_5 + 4 \xi_4 \xi_5 - 2 \xi_2 \xi_6 - 8 \xi_4 \xi_6 + 6 \xi_5 \xi_6 = 0.
\end{dmath}
This is the realization of Theorem~\ref{th:criallel} in the case of 7-by-7 matrices. 

However, further simplification is possible. Observe that the difference between the second and the third equation in \eqref{con7} is nothing but \[ 4(\xi_1+\xi_2-\xi_5-\xi_6)(\xi_2+\xi_3-\xi_4-\xi_5).\]  

So, for $C(A)$ to consist of the concentric ellipses it is necessary that either $\xi_1+\xi_2-\xi_5-\xi_6=0$ or $\xi_2+\xi_3-\xi_4-\xi_5=0$. (Note that the same conclusion also follows from Corollary~\ref{th:midel7}.) 

The former condition, when combined with \eqref{con7}, upon applying the Gr\"obner basis transformation, simplifies to 
\begin{equation}\label{con71}
\begin{aligned}
&\xi_2^2 - \xi_3^2 - \xi_4^2 - \xi_5^2 - 2 \xi_6^2 - 6 \xi_3 \xi_2 + 2 \xi_4 \xi_2 + 2 \xi_5 \xi_2 - 2 \xi_3 \xi_4 \\ 
&+ 6 \xi_3 \xi_5 - 2 \xi_4 \xi_5 + 4 \xi_3 \xi_6 + 4 \xi_4 \xi_6 - 4 \xi_5 \xi_6 = 0, \\
&\xi_1 + \xi_2 - \xi_5 - \xi_6 = 0.
\end{aligned}
\end{equation}

So, the following result holds. 
\begin{cor}\label{th:conel7}Let $A\in\Rec_7$. Then $C(A)$ consists of three concentric ellipses and the point at the origin if and only if either \eqref{con71} holds, or $\xi_2+\xi_3=\xi_4+\xi_5$ and the first two equation in \eqref{con7} are satisfied.  \end{cor} 

Note that for $\xi_j$ as in \eqref{ex71} the first equality in \eqref{con71} fails while apparently $\xi_2+\xi_3\neq\xi_4+\xi_5$. So, for the curve in Fig.~\ref{fig:1} the interior and exterior components thereof are not elliptical (despite looking very close to such).  

To provide an example of a matrix with the Kippenhahn curve consisting of three concentric ellipses, let 
\[\xi_1 = \xi_2 = \xi_5 = \xi_6 = 1,\ \xi_3 = 2, \ \xi_4 = 0\]

\begin{figure}[H]
    \centering
    \includegraphics[width=0.5\textwidth]{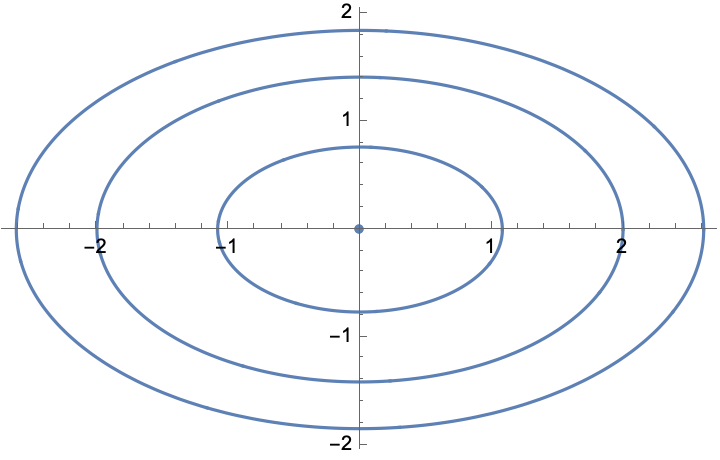}
\caption{$\vec{\xi} = \left\{1, 1, 2,0, 1, 1\right\}$}
\label{fig:conel7}
\end{figure}

This choice fits conditions of Theorem~\ref{th:crionel} for the case of 7-by-7 matrices and all values of $X=X_j$, $j=1,2,3$, as well as Corollary~\ref{th:conel7}. The respective Kippenhahn curve is plotted in Fig~\ref{fig:conel7}. 

\section{Shifted Ellipses. Criteria for $n=7$} \label{s:shif7}
For $n=4$ condition \eqref{pash} holds if and only if $C(A)$ actually is the union of two shifted ellipses. The respective criterion was obtained in \cite[Theorem 6.1]{GeS2}. In the case $n=5$, \eqref{pash} is equivalent to $C(A)=E\cup (-E)\cup\{0\}$, and the criterion for that to happen is \cite[Theorem 6]{JiangS}. The same paper contains a detailed treatment of the case $n=6$ in which  \eqref{pash} implies 
\[ C(A)=E\cup (-E)\cup E_0, \] where $E_0$ 
also is an ellipse, but centered at the origin \cite[Theorem 9]{JiangS}. It can be verified that Theorem~\ref{th:shifel} provides a unified approach to all these cases. When applied to the $n=7$ setting, equating the defined by \eqref{roe} polynomials $R_o,R_e$ with zero, as described in Section~\ref{s:shifel}, leads to the following criterion.
\begin{prop}\label{th:crishift7}Let $A\in\Rec_7$. Then \eqref{pash} holds if and only if for some $p,X,C>0$:  
\eq{pX7} p\pm X\in\left\{ \pm\sqrt{2+\sqrt{2}},\pm\sqrt{2},\pm\sqrt{2-\sqrt{2}}\right\}\en 
and \end{prop}  \begin{dmath} \label{crish}
2 C \left(15 p^4+6 p^2 \left(5 X^2-6\right)+3 X^4-12 X^2+10\right)-2 \xi _2-4 \xi _3-4 \xi _4-2 \xi _5-6 \xi _6-2 \xi _1 \left(p^4+p^2 \left(6 X^2-4\right)+X^4-4 X^2+3\right)-2 \left(\xi _2+\xi _3+\xi _4+\xi _5\right) \left(p^4+p^2 \left(6 X^2-3\right)+X^4-3 X^2\right)-2 \xi _6 \left(p^4+p^2 \left(6 X^2-4\right)+X^4-4 X^2\right)=0,\\
2 C \left(5 p^2+3 X^2-6\right)-2\left(p^2+X^2-2\right)(\xi_1+\xi_6)+(3-2(p^2+X^2)) \left(\xi _2+\xi _3+\xi _4+\xi _5\right)\nolinebreak =0,\\ 
3C^2(5p^2+X^2-2)-C(6p^2+2X^2-3)(\xi_2+\xi_3+\xi_4+\xi_5)-2C(3p^2+X^2-2)(\xi_1+\xi_6)+(p^2+X^2-1)(\xi_1\xi_5+\xi_2\xi_4+\xi_2\xi_6+\xi_3\xi_5+\xi_3\xi_6)+(p^2+X^2-2)(\xi_1\xi_3+\xi_1\xi_4+\xi_1\xi_6+\xi_4\xi_6) \nolinebreak =0,\\
3 C^2- 2C(\xi_1+\xi_2+\xi_3+\xi_4+\xi_5+\xi_6)+\xi_1(\xi_3+\xi_4+\xi_5+\xi_6)+\xi_2(\xi_4+\xi_5+\xi_6) +\xi_3(\xi_5+\xi_6)+\xi_4\xi_6 \nolinebreak =0,\\
C^3-C^2(\xi_1+\xi_2+\xi_3+\xi_4+\xi_5+\xi_6)+C(\xi_1(\xi_3+\xi_4+\xi_5+\xi_6)+\xi_2(\xi_4+\xi_5+\xi_6)+
\xi_3(\xi_5+\xi_6)+\xi_4\xi_6)-(\xi_1\xi_3(\xi_5+\xi_6)+(\xi_1+\xi_2)\xi_4\xi_6) \nolinebreak =0.
\end{dmath}
{\em If these conditions hold, then in \eqref{pash} $p\pm X$ and $-p\pm X$ are the foci of $E$ and $-E$, respectively, and $C$ is the square of the length of their minor half-axes.} 

It was observed in \cite{JiangS} that for $n=4,5,6$ the inclusion \eqref{pash} can materialize only for ellipses $E$ having foci of opposite sign. It was conjectured there that this property persists for $n>6$. Based on Proposition~\ref{th:crishift7} we can show that it is still the case for $n=7$.
\begin{thm}\label{th:noshift7} Let $A\in\Rec_7$ be such that \eqref{pash} holds. Then the foci of $E$ (and thus of $-E$ as well) are of the opposite sign.\end{thm}
\begin{proof} Excluding $C$ from the first two equations in \eqref{crish}, we end up with 
\begin{dmath}\label{linxi} 
2\left(10 p^6+12 p^4 X^2-40 p^4+10 p^2 X^4-40 p^2 X^2+43 p^2+X^2-2\right)(\xi _1+\xi _6) + \left(20 p^6+24 p^4 X^2-75 p^4+20 p^2 X^4-66 p^2 X^2+82 p^2-3 X^4+14 X^2-18\right)(\xi _2+\xi _4) + \left(20 p^6+24 p^4 X^2-75 p^4+20 p^2 X^4-66 p^2 X^2+72 p^2-3 X^4+8 X^2-6\right)(\xi _3+\xi _5) = 0.
\end{dmath} 
Direct computations show that for all choices of $X,p$ in \eqref{pX7} with $p>X$ the coefficients in \eqref{linxi} are positive. 

Namely, for \[ p=\frac{1}{2} \left(\sqrt{\sqrt{2}+2}+\sqrt{2}\right), \quad X=\frac{1}{2} \left(\sqrt{\sqrt{2}+2}-\sqrt{2}\right), \] \eqref{linxi} takes the form 
\begin{multline*}
    \left(2 \left(\sqrt{2}+3\right)+3 \sqrt{2 \left(\sqrt{2}+2\right)}\right) (\xi _1+\xi_6)\\ +\left(4 \sqrt{2}+\sqrt{46 \sqrt{2}+68}+6\right) (\xi _2+\xi_5)\\+2 \left(\sqrt{2}+\sqrt{7 \sqrt{2}+10}+2\right) (\xi _3+\xi_4)= 0.
\end{multline*}

For \[ p=\frac{1}{2} \left(\sqrt{2-\sqrt{2}}+\sqrt{\sqrt{2}+2}\right), \quad X=\frac{1}{2} \left(\sqrt{\sqrt{2}+2}-\sqrt{2-\sqrt{2}}\right), \] \eqref{linxi} takes the form 
\begin{multline*}
    \left(\sqrt{2}+2\right)(\xi _1+\xi_6)+\sqrt{2} (\xi _2+\xi_5)+2\left(\sqrt{2}+1\right) (\xi _3+\xi_4)= 0.
\end{multline*}
Finally, for \[ p=\frac{1}{2} \left(\sqrt{2}+\sqrt{2-\sqrt{2}}\right), \quad X=\frac{1}{2} \left(\sqrt{2}-\sqrt{2-\sqrt{2}}\right), \] \eqref{linxi} takes the form 
\begin{multline*}
    \left(6-2\sqrt{2}+3 \sqrt{4-2 \sqrt{2}}\right) (\xi _1+\xi_6)\\ +\left(6-4 \sqrt{2}-2 \sqrt{2-\sqrt{2}}+3 \sqrt{4-2 \sqrt{2}}\right) (\xi _2+\xi_5)\\+2 \left(2-\sqrt{2}-\sqrt{2-\sqrt{2}}+\sqrt{4-2 \sqrt{2}}\right)(\xi _3+\xi_4)= 0.
\end{multline*}

Therefore, equation \eqref{linxi} does not admit non-trivial non-negative solutions. The same therefore holds for the system \eqref{crish}. \end{proof} 

According to Theorem~\ref{th:noshift7}, we are left with the three cases when \eqref{pX7} is subjected to the additional requirement $p<X$.

Equation \eqref{linxi} still holds, but (as another run of numerical computations shows) its coefficients now vary in signs. So, the reasoning from the proof of Theorem~\ref{th:noshift7} does not apply. 

Instead, completing the procedure described in Section~\ref{s:shifel}, for each admissible pair $\{p,X\}$ we can solve \eqref{linxi} for $C$ and plug the solution into \eqref{crish}. This yields the criteria for \eqref{pash} to hold consisting of four homogeneous polynomial equations in $\xi_1,\ldots,\xi_6$: one linear, two quadratic, and one cubic. Our efforts to generate exact numerical examples by solving these systems with the use of \textit{Mathematica} were unsuccessful.

However, due to a still relatively low value of $n$, an alternative approach is available. 

\begin{thm} \label{th:altcri} For $A\in \Rec_7$, $C(A)$ contains a pair of shifted ellipses if and only if 
the polynomial \eqref{P7} factors as
\eq{p7shif} \left(\zeta -\left(C_0+\rho  X_0^2\right)\right) \left(\zeta ^2-2 \zeta  \left(C+\rho  \left(p^2+X^2\right)\right)+\left(C+\rho  \left(X^2-p^2\right)\right)^2\right) \en 
with  $p<X$, $p\pm X$ satisfying \eqref{pX7}, $X_0$ being the positive eigenvalue of $A$ different from $X\pm p$, and either \newline \em{(i)}  $\xi_2\xi_3$ = 0,
\begin{align}\label{eq:quad}
         \left(\xi _1+\xi _2\right){}^2-\left(\xi _3+\xi _4+\xi _5+\xi _6\right) \left(\xi _1+\xi _2\right)+\xi _3 \xi _5+\xi _3 \xi _6+\xi _4 \xi _6 = 0,
     \end{align}
     \begin{align}\label{eq:c0c}
        C = \xi _1+\xi _2, \quad C_0 = \xi _3+\xi _4+\xi _5+\xi _6-\left(\xi _1+\xi _2\right),
     \end{align}
or (ii)  $\xi_4\xi_5$ = 0, 
\begin{align}\label{eq:quadi}
\left(\xi _5+\xi _6\right){}^2-\left(\xi _1+\xi _2+\xi _3+\xi _4\right) \left(\xi _5+\xi _6\right)+\xi _1 \xi _3+\xi_1 \xi _4+\xi _2 \xi _4 = 0,
\end{align}
\begin{align}\label{eq:c0ci}
        C = \xi _5+\xi _6, \quad C_0 = \xi _1+\xi _2+\xi _3+\xi_4-\left(\xi _5+\xi _6\right).
     \end{align}
If these conditions hold, then in fact \eq{ee0}C(A)=E\cup(-E)\cup E_0\cup\{0\},\en  where $E$ is the ellipse with the foci $p\pm X$ and the minor half-axis $\sqrt{C}$, and $E_0$ is the ellipse with the foci $\pm X_0$ and the minor half-axis $\sqrt{C_0}$.  
\end{thm}
\begin{proof}
As was shown in Section~\ref{s:shifel}, for $A\in\Rec_n$ the Kippenhahn curve $C(A)$ contains a pair of shifted ellipses $E\cup(-E)$ if and only if the polynomial $P_n(\zeta)$ is divisible by the quadratic \eqref{shell} which is the second factor in the right-hand side of \eqref{p7shif}. Since $P_7$ has degree three in $\zeta$, the quotient of its division by \eqref{shell} is linear. The respective component of $C(A)$ therefore is also an ellipse, say $E_0$, and its foci are the non-zero eigenvalues of $A$ different from $\pm p\pm X$. These are exactly $\pm X_0$, implying that $E_0$ is centered at the origin. According to \eqref{eq0}, the factorization of $P_7$ has the form \eqref{p7shif}. 

Condition \eqref{p7shif} with $p,X,X_0$ as described and $C,C_0$ considered as unknowns along with $\xi_1,\ldots,\xi_6$ is already a criterion for $C(A)$ to contain a pair of shifted ellipses. To derive \eqref{eq:c0c} and \eqref{eq:quad}, observe that \eqref{pash} implies that the joint tangent lines of $E$ and $-E$ are multiple tangent lines of $C(A)$. According to \cite[Proposition 4]{JiangS}, for $A\in\Rec_n$ this is only possible if $\xi_k=0$ for some $k\in\{2,\ldots,n-2\}$ and the spectra of the left-upper $k$-by-$k$ block $B_1$ and the right lower $(n-k)$-by-$(n-k)$ block $B_2$ of $\im A$ overlap. 
Moreover, the ordinates of the tangent lines in question are nothing but the points in $\sigma(B_1)\cap\sigma(B_2)$.

For $n=7$ this yields four cases, in each of which the characteristic polynomials $P_1,P_2$ of $B_1,B_2$ can be computed with the use of Proposition~\ref{th:exptri}. 

{\sl Case 1.} $\xi_2=0$.  Then $P_1(\lambda)=\zeta-\xi_1$, \[ P_2(\lambda)=\lambda(\zeta^2-(\xi_3+\xi_4+\xi_5+\xi_6)\zeta+\xi_3\xi_5+\xi_3\xi_6+\xi_4\xi_6),\] where $\zeta=\lambda^2$ as before. Consequently, $P_1$ and $P_2$ have common roots  if and only if $P_2$ vanishes under the substitution $\zeta=\xi_1$:
\[  \xi_1^2-(\xi_3+\xi_4+\xi_5+\xi_6)\xi_1+(\xi_3\xi_5+\xi_3\xi_6+\xi_4\xi_6)=0, \]
and these roots are $\pm\sqrt{\xi_1}$. This agrees with \eqref{eq:quad} and the expression for $C$ in \eqref{eq:c0c} considering that $\xi_2=0$. The remaining non-zero roots of $P_2$ are $\pm\sqrt{\xi_3+\xi_4+\xi_5+\xi_6-\xi_1}$ which agrees with the second formula in 
 \eqref{eq:c0c}.

{\sl Case 2.} $\xi_3=0$. Then 
\[ P_1(\lambda) = \lambda(\xi_1+\xi_2-\zeta), \quad
P_2(\lambda) = \zeta^2-\zeta(\xi_4+\xi_5+\xi_6)+\xi_4\xi_6, \]
and $P_1,P_2$ have common roots if and only if $P_2$ vanishes under the substitution $\zeta=\xi_1+\xi_2$. But this is exactly \eqref{eq:quad} when $\xi_3=0$. Formulas \eqref{eq:c0c} follow in a similar manner. 

This takes care of (i). The remaining two cases corresponding to (ii) can be treated in the same manner. Alternatively, it suffices to observe that a transpositional similarity yields the substitutions \eq{subst}  \xi_1\longleftrightarrow\xi_6, \quad \xi_2\longleftrightarrow\xi_5, \quad \xi_3\longleftrightarrow\xi_4 \en 
under which \eqref{eq:quad}--\eqref{eq:c0c} go into \eqref{eq:quadi}--\eqref{eq:c0ci}. \end{proof} 

\section{Shifted Ellipses. Explicit descriptions for $n=7$}\label{s:shif7ex}
There are three possible configurations falling under the setting of Theorem~\ref{th:altcri}, corresponding to the three choices of $X_0$. The next statements deal with these three possibilities, one at a time.
\begin{thm}\label{th:7central}For $A\in \Rec_7$, $C(A)$ contains an ellipse $E$ with the foci $\sqrt{2+\sqrt{2}}$ and $-\sqrt{2-\sqrt{2}}$ (equivalently, $\sqrt{2-\sqrt{2}}$ and $-\sqrt{2+\sqrt{2}}$) if and only if the 6-tuple $(\xi_1,\ldots,\xi_6)$ for some $t>0$ equals $t{\vec{\xi}}_1$ or $t{\vec{\xi}}_2$ with \eq{Xi} {\vec{\xi}}_1=(2(\sqrt{2}-1), 3-2\sqrt{2}, 0, 1,0,1), \
{\vec{\xi}}_2=(\sqrt{2}+1, 0, \sqrt{2}+1, 0, \sqrt{2}-1, 2).\en  
 
If this is the case, then in fact $C(A)$ is given by \eqref{ee0}, with $E_0$ being an ellipse with the foci $\pm\sqrt{2}$. Moreover, the minor axes of $\pm E$ and $E_0$ are of the same length. 
\end{thm} 
\begin{proof}The configuration in question corresponds to
\eq{midx0} X_0=\sqrt{2}, \quad X\pm p=\sqrt{2\pm\sqrt{2}}.\en
Plugging \eqref{midx0} into \eqref{p7shif} and equating its coefficients with the respective ones in \eqref{P7} yields 

\begin{align*}
&C^2 C_0 - \xi_2 \xi_4 \xi_6 - \xi_1 (\xi_4 \xi_6 + \xi_3 (\xi_5 + \xi_6)) = 0, \\
&2 C (C + \sqrt{2} C_0) - \xi_2 \xi_4 - \xi_3 \xi_5 \\&- (\xi_2 + \xi_3 + 2 \xi_4) \xi_6 - \xi_1 (2 \xi_3 + \xi_4 + \xi_5 + 2 \xi_6) = 0, \\
&4 \sqrt{2} C + 2 C_0 - 3 \xi_1 - \xi_2 - 2 \xi_3 - 2 \xi_4 - \xi_5 - 3 \xi_6 = 0, \\
&-C (C + 2 C_0) + \xi_3 \xi_5 + (\xi_3 + \xi_4) \xi_6 + \xi_2 (\xi_4 + \xi_5 + \xi_6) \\&+ \xi_1 (\xi_3 + \xi_4 + \xi_5 + \xi_6) = 0, \\
&-2 (\sqrt{2} + 2) C - 4 C_0 + 4 \xi_1 + 3 \xi_2 + 3 \xi_3 + 3 \xi_4 + 3 \xi_5 + 4 \xi_6 = 0, \\
&2 C + C_0 - \xi_1 - \xi_2 - \xi_3 - \xi_4 - \xi_5 - \xi_6 = 0.
\end{align*}
Using \eqref{eq:c0c} to replace $C,C_0$ by their expressions in terms of $\xi_j$ and augmenting this system with \eqref{eq:quad}, we end up with 

\begin{align}\label{sys}
&-\xi _1^3 - 3 \xi _2 \xi _1^2 + \xi _3 \xi _1^2 + \xi _4 \xi _1^2 + \xi _5 \xi _1^2 + \xi _6 \xi _1^2 - 3 \xi _2^2 \xi _1 + 2 \xi _2 \xi _3 \xi _1 \nonumber \\
&\quad + 2 \xi _2 \xi _4 \xi _1 + 2 \xi _2 \xi _5 \xi _1 - \xi _3 \xi _5 \xi _1 + 2 \xi _2 \xi _6 \xi _1 - \xi _3 \xi _6 \xi _1 - \xi _4 \xi _6 \xi _1 \nonumber \\
&\quad - \xi _2^3 + \xi _2^2 \xi _3 + \xi _2^2 \xi _4 + \xi _2^2 \xi _5 + \xi _2^2 \xi _6 - \xi _2 \xi _4 \xi _6 = 0, \nonumber \\
&-2 \sqrt{2} \xi _1^2 + 2 \xi _1^2 - 4 \sqrt{2} \xi _2 \xi _1 + 4 \xi _2 \xi _1 + 2 \sqrt{2} \xi _3 \xi _1 - 2 \xi _3 \xi _1 \nonumber \\
&\quad + 2 \sqrt{2} \xi _4 \xi _1 - \xi _4 \xi _1 + 2 \sqrt{2} \xi _5 \xi _1 - \xi _5 \xi _1 + 2 \sqrt{2} \xi _6 \xi _1 - 2 \xi _6 \xi _1 \nonumber \\
&\quad - 2 \sqrt{2} \xi _2^2 + 2 \xi _2^2 + 2 \sqrt{2} \xi _2 \xi _3 + 2 \sqrt{2} \xi _2 \xi _4 - \xi _2 \xi _4 + 2 \sqrt{2} \xi _2 \xi _5 \nonumber \\
&\quad - \xi _3 \xi _5 + 2 \sqrt{2} \xi _2 \xi _6 - \xi _2 \xi _6 - \xi _3 \xi _6 - 2 \xi _4 \xi _6 = 0, \nonumber \\
&4 \sqrt{2} \xi _1 - 5 \xi _1 + 4 \sqrt{2} \xi _2 - 3 \xi _2 + \xi _5 - \xi _6 = 0, \nonumber \\
&\xi _1^2 + 2 \xi _2 \xi _1 - \xi _3 \xi _1 - \xi _4 \xi _1 - \xi _5 \xi _1 - \xi _6 \xi _1 + \xi _2^2 - 2 \xi _2 \xi _3 \nonumber \\
&\quad - \xi _2 \xi _4 - \xi _2 \xi _5 + \xi _3 \xi _5 - \xi _2 \xi _6 + \xi _3 \xi _6 + \xi _4 \xi _6 = 0, \nonumber \\
&-2 \sqrt{2} \xi _1 + 4 \xi _1 - 2 \sqrt{2} \xi _2 + 3 \xi _2 - \xi _3 - \xi _4 - \xi _5 = 0, \nonumber \\
&(\xi _1 + \xi _2)^2 - (\xi _3 + \xi _4 + \xi _5 + \xi _6)(\xi _1 + \xi _2) \nonumber \\
&\quad + \xi _3 \xi _5 + \xi _3 \xi _6 + \xi _4 \xi _6 = 0
\end{align}
Following the logic of Theorem~\ref{th:altcri}, we now have to treat the cases when one of the variables $\xi_2,\ldots,\xi_5$ is equal to zero. 

Setting $\xi_2=0$ and applying Mathematica Gr\"obner basis procedure allows to rewrite \eqref{sys} in an equivalent form 
\begin{align}\label{gb2}
    &34 \xi _5^3+9 \sqrt{2} \xi _6 \xi _5^2-31 \xi _6 \xi _5^2+65 \sqrt{2} \xi _6^2 \xi _5-86 \xi _6^2 \xi _5+32 \sqrt{2} \xi _6^3-46 \xi _6^3=0, \nonumber \\&104 \sqrt{2} \xi _5^2+158 \xi _5^2-5 \sqrt{2} \xi _6 \xi _5-\xi _6 \xi _5-\sqrt{2} \xi _6^2-10 \xi _6^2+49 \xi _4 \xi _6=0, \nonumber \\&48 \sqrt{2} \xi _5^2+88 \xi _5^2+49 \xi _4 \xi _5-40 \sqrt{2} \xi _6 \xi _5-8 \xi _6 \xi _5-8 \sqrt{2} \xi _6^2+18 \xi _6^2=0, \nonumber \\& 7 \xi _3+7 \xi _4+6 \sqrt{2} \xi _5+11 \xi _5-6 \sqrt{2} \xi _6-4 \xi _6=0, \nonumber \\& 7 \xi _1+4 \sqrt{2} \xi _5+5 \xi _5-4 \sqrt{2} \xi _6-5 \xi _6 = 0.
\end{align}
All positive solutions of the system \eqref{gb2} are indeed of the form $t{\vec{\xi}}_2$ with $t>0$. 

Observe also that, according to \eqref{eq:c0c}, in this case $C=C_0=(\sqrt{2}+1)t$. 

Similarly, letting $\xi_3=0$ in \eqref{sys} with the use of Gr\"obner's procedure we can rewrite it as 
\begin{align}\label{gb1}
&\xi_5 (\xi_5 + \xi_6) (\xi_5^2 - 3 \sqrt{2} \xi_6 \xi_5 - 4 \xi_6 \xi_5 - 2 \sqrt{2} \xi_6^2 + 2 \xi_6^2) = 0, \nonumber \\
&204 \sqrt{2} \xi_5^3 + 123 \xi_5^3 - 903 \sqrt{2} \xi_6 \xi_5^2 - 1575 \xi_6 \xi_5^2 - 1107 \sqrt{2} \xi_6^2 \xi_5 - 2356 \xi_6^2 \xi_5 \nonumber \\
&\quad - 658 \xi_6^3 + 658 \xi_4 \xi_6^2 = 0, \nonumber \\
&6 \sqrt{2} \xi_5^2 + 16 \xi_5^2 + 23 \xi_4 \xi_5 + 4 \sqrt{2} \xi_6 \xi_5 + 26 \xi_6 \xi_5 - 2 \sqrt{2} \xi_6^2 + 10 \xi_6^2 \nonumber \\
&\quad + 2 \sqrt{2} \xi_4 \xi_6 - 10 \xi_4 \xi_6 = 0, \nonumber \\
&23 \xi_4^2 - 5 \sqrt{2} \xi_6 \xi_4 - 21 \xi_6 \xi_4 - 15 \sqrt{2} \xi_5^2 - 17 \xi_5^2 + 5 \sqrt{2} \xi_6^2 - 2 \xi_6^2 \nonumber \\
&\quad - 10 \sqrt{2} \xi_5 \xi_6 - 19 \xi_5 \xi_6 = 0, \nonumber \\
&3 \xi_2 + 4 \sqrt{2} \xi_4 - 5 \xi_4 + 2 \sqrt{2} \xi_5 - \xi_5 + 2 \sqrt{2} \xi_6 - 4 \xi_6 = 0, \nonumber \\
&3 \xi_1 - 4 \sqrt{2} \xi_4 + 3 \xi_4 - 2 \sqrt{2} \xi_5 - 2 \sqrt{2} \xi_6 + 3 \xi_6 = 0.
\end{align}
The general positive solution of the latter is $t{\vec{\xi}}_1,\ t>0$, with $t$ actually equal to both $C$ and $C_0$.  

There is no need to consider $\xi_4=0$ and $\xi_5=0$ separately, since substitutions \eqref{subst} reduce
these cases to already treated ones. Moreover, no new solutions emerge: $\xi_4=0$ yields the same solution as $\xi_2=0$, and $\xi_5=0$ --- the same as $\xi_3=0$. The proof is thus complete. 
\end{proof} 

To illustrate, $C(A)$ corresponding to $t=1$ in \eqref{Xi} is plotted in Figures~3 and 4, respectively. 
\begin{figure}[H]\label{c1} 
    \centering
    \includegraphics[width=0.5\linewidth]{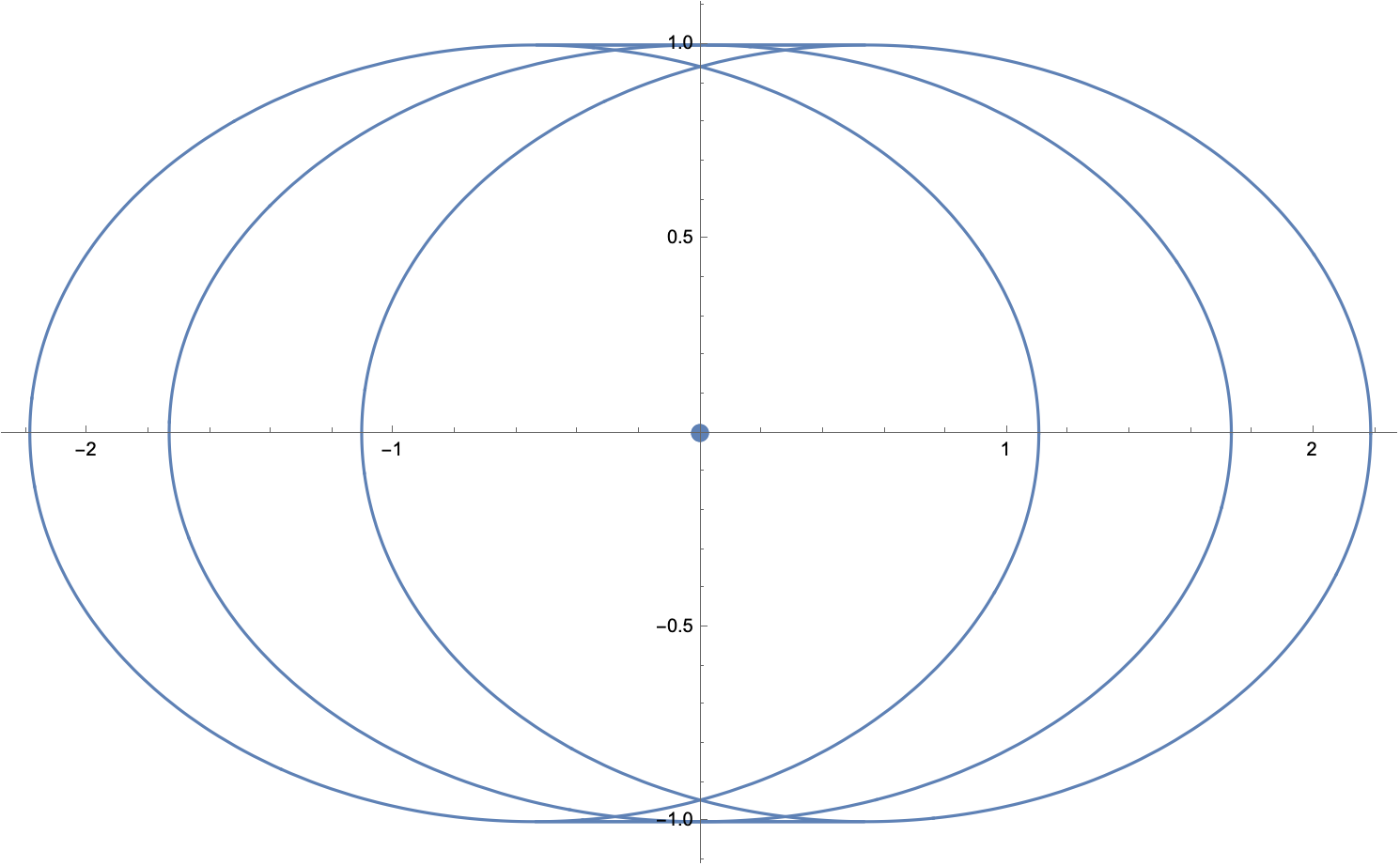}
    \caption{$\vec{\xi}\to\{ 2 \left(\sqrt{2}-1\right), 3-2 \sqrt{2}, 0, 1, 0, 1\}$}
\end{figure}
 \begin{figure}[H]\label{c2} 
    \centering
    \includegraphics[width=0.5\linewidth]{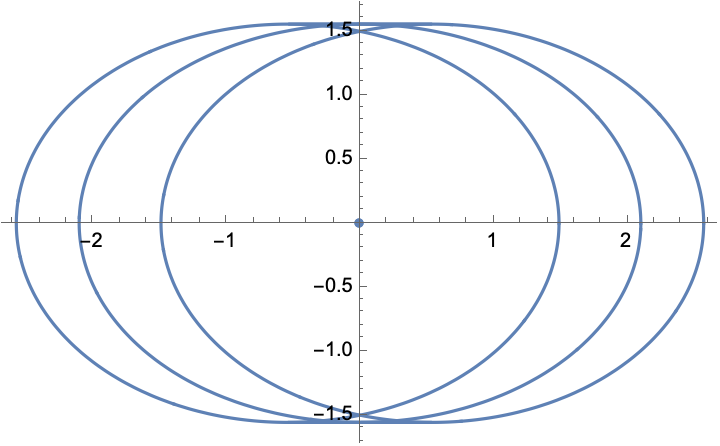}
    \caption{$\vec{\xi}\to\left\{\sqrt{2}+1, 0, \sqrt{2}+1, 0, \sqrt{2}-1, 2\right\}$}
\end{figure}
\begin{thm}\label{th:7inner}
For $A\in \Rec_7$, $C(A)$ contains an ellipse $E$ with the foci $\sqrt{2+\sqrt{2}}$ and $-\sqrt{2}$ (equivalently, $\sqrt{2}$ and $-\sqrt{2+\sqrt{2}}$) if and only if the 6-tuple $(\xi_1,\ldots,\xi_6)$ for some $t>0$ and $j=1,\ldots,8$ equals $t\vec{\xi}_j$ with 

\begin{equation}\label{Ximid}
\begin{aligned}
\vec{\xi}_1 &= \{1,0,-7 \sqrt{2}+4 \sqrt{10-7 \sqrt{2}}+9,2 (\sqrt{2}-2 \sqrt{2 (\sqrt{2}+2)}+4), \\
&\quad -8 \sqrt{2}+6 \sqrt{4-2 \sqrt{2}}+5,2 (2 \sqrt{2}+\sqrt{20-14 \sqrt{2}}-3)\} \\
&\approx \{1, 0, 0.368606, 0.375923, 0.180645, 0.553537\}, \\
\vec{\xi}_2 &= \{1,0,-2 \sqrt{\sqrt{2}+2}+\frac{3}{\sqrt{2}}+2,2 \sqrt{2-\sqrt{2}}+\frac{1}{\sqrt{2}}-2, \\
&\quad 2 (\sqrt{2}-2 \sqrt{2 (\sqrt{2}+2)}+4),-6 \sqrt{2}+4 \sqrt{7 \sqrt{2}+10}-9\} \\
&\approx \{1, 0, 0.425802, 0.237841, 0.375923, 0.358259\}, \\
\vec{\xi}_3 &\approx \{1,1.69724,0,1.01396,1.19026,0.790003\}, \\
\vec{\xi}_4 &\approx \{1,1.51367,0,1.0412,0.944947,0.900544\}.
\end{aligned}
\end{equation}
and, for $i=1,\ldots 4$, $\vec{\xi}_{i+4}$ obtained from $\vec{\xi}_i$ via the permutation
\eq{perm} 
\xi_j \longleftrightarrow \xi_{7-j} \quad (j = 1,\ldots,6).
\en

If this is the case, then in fact $C(A)$ is given by \eqref{ee0}, with $E_0$ being an ellipse with the foci $\pm\sqrt{2-\sqrt{2}}$. Moreover, the minor half-axes of $\pm E$ and $E_0$ have length 
\[
\left\{1,\sqrt{15-9 \sqrt{2}-4 \sqrt{20-14 \sqrt{2}}}\right\}t, \quad \left\{1,\sqrt{2 (\sqrt{2 (\sqrt{2}+2)}-\sqrt{2}-1)}\right\}t,
\]
\[
\left\{1.64233,0.544959\right\}t, \quad \left\{1.58546,0.610752\right\}t,
\]
respectively.
\end{thm} 
\begin{proof} This configuration corresponds to $X+p=\sqrt{2+\sqrt{2}},\ X-p=\sqrt{2}$, and $X_0=\sqrt{2-\sqrt{2}}$. For this choice of the parameters, case (i) of the criterion of Theorem~\ref{th:altcri} for $\xi_2=0$ takes the form
\begin{align*}
&\xi_1 \left(-\xi_1^2 + (\xi_3 + \xi_4 + \xi_5 + \xi_6) \xi_1 - \xi_4 \xi_6 - \xi_3 (\xi_5 + \xi_6)\right) = 0, \\
&-\left((\sqrt{2} + 2 \sqrt{2 (\sqrt{2} + 2)} - 2) \xi_1^2\right) \\
&\quad + \left(2 (\sqrt{2 (\sqrt{2} + 2)} - 1) \xi_3 - \xi_4 - \xi_5 - 2 \xi_6 + 2 \sqrt{2 (\sqrt{2} + 2)} (\xi_4 + \xi_5 + \xi_6)\right) \xi_1 \\
&\quad - 2 \xi_4 \xi_6 - \xi_3 (\xi_5 + \xi_6) = 0, \\
&\left(-2 \sqrt{2} + 4 \sqrt{2 - \sqrt{2}} - 7\right) \xi_1 + 2 (\sqrt{2} + 1) \xi_3 + 2 (\sqrt{2} + 1) \xi_4 \\
&\quad + (2 \sqrt{2} + 3) \xi_5 + (2 \sqrt{2} + 1) \xi_6 = 0, \\
&\xi_1^2 - (\xi_3 + \xi_4 + \xi_5 + \xi_6) \xi_1 + \xi_4 \xi_6 + \xi_3 (\xi_5 + \xi_6) = 0, \\
&\left(3 \sqrt{2} - 2\sqrt{2(2+\sqrt{2}} + 4\right) \xi_1 + \frac{\sqrt{2} \xi_3 + \sqrt{2} \xi_4 + \sqrt{2} \xi_5 + 2 (\sqrt{2} - 1) \xi_6}{\sqrt{2} - 2} = 0.
\end{align*}
If $\xi_1=0$ then the first equation is satisfied but the last one implies $\xi_j=0$ for $j=3,\ldots,6$ which is a trivial (and thus unacceptable) solution. So, the factor $\xi_1$ in the first equation can be dropped, with the resulting equation exactly the same as the fourth one. This leaves us with two linear and two quadratic equations in five unknowns.  Solving the linear equations for $\xi_5,\xi_6$, we obtain:
\eq{xi56}
\begin{aligned}
    \xi_5 &= \left(4 \left(\sqrt{2}+3\right)-6 \sqrt{2 \left(\sqrt{2}+2\right)}\right) \xi _1-\left(\sqrt{2}+1\right) \left(\xi _3+\xi _4\right),\\
    \xi_6 &= \left(-8 \sqrt{2}+2 \sqrt{46 \sqrt{2}+68}-13\right) \xi _1+\left(\sqrt{2}+1\right) \left(\xi _3+\xi _4\right).
\end{aligned}
\en 
Plugging \eqref{xi56} into the remaining equations and setting temporarily $\xi_1=1$, we end up with the system of two quadratic equations in two variables $\xi_3,\xi_4$. This system has exactly two positive solutions, namely 
\[ \xi_3=-7 \sqrt{2}+4 \sqrt{10-7 \sqrt{2}}+9,\ \xi_4= 2 \left(\sqrt{2}-2 \sqrt{2 \left(\sqrt{2}+2\right)}+4\right) \]
and 
\[ \xi_3=-2 \sqrt{\sqrt{2}+2}+\frac{3}{\sqrt{2}}+2,\ \xi_4=2 \sqrt{2-\sqrt{2}}+\frac{1}{\sqrt{2}}-2. \]
Recovering $\xi_5,\xi_6$ from \eqref{xi56}, we arrive at ${\vec{\xi}}_1$ and ${\vec{\xi}}_2$ from \eqref{Ximid}. 

A similar reasoning in the case (i) when $\xi_3=0$ also generates two one-parameter families of the solutions, corresponding to ${\vec{\xi}}_3, {\vec{\xi}}_4$ in \eqref{Ximid}. Note that only an approximate solution is available here. 

Finally, case (ii) of the criterion from Theorem~\ref{th:altcri} follows immediately with the use of \eqref{perm}.
\end{proof} 

The Kippenhahn curve of the matrix $A$ corresponding to $\xi_j$ given by $\vec{\xi}_1$ and $t = 1$ is shown below.
\begin{figure}[H]
    \centering
    \includegraphics[width=0.5\linewidth]{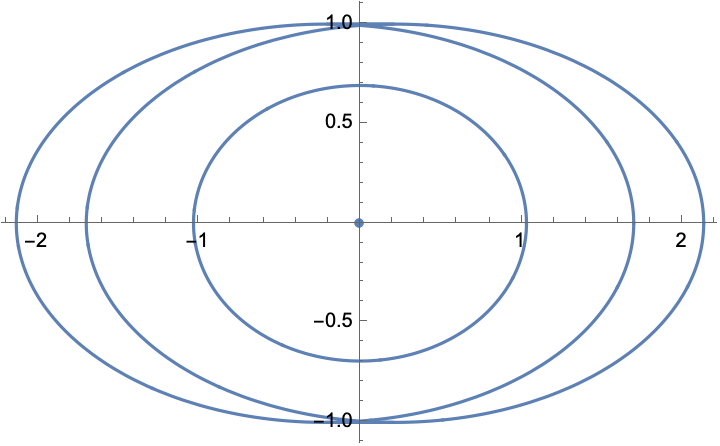}
    \begin{minipage}{1.0\textwidth}
        \centering
        \caption{
            $\vec{\xi} \to \{1, 0, -7 \sqrt{2} + 4 \sqrt{10 - 7 \sqrt{2}} + 9, 2 (\sqrt{2} - 2 \sqrt{2 (\sqrt{2} + 2)} + 4), -8 \sqrt{2} + 6 \sqrt{4 - 2 \sqrt{2}} + 5, 2 (2 \sqrt{2} + \sqrt{20 - 14 \sqrt{2}} - 3)\}$
        }\end{minipage}
\end{figure}

\begin{thm}\label{th:7outer}For $A\in \Rec_7$, $C(A)$ contains an ellipse $E$ with the foci $\sqrt{2-\sqrt{2}}$ and $-\sqrt{2}$ (equivalently, $\sqrt{2}$ and $-\sqrt{2-\sqrt{2}}$) if and only if the 6-tuple $(\xi_1,\ldots,\xi_6)$ for some $t>0$ and $j=1,\ldots,6$ equals $t{\vec{\xi}}_j$ with 
\[ 
\begin{aligned}
{\vec{\xi}}_1 &= \{1,0,7 \sqrt{2}-4 \sqrt{7 \sqrt{2}+10}+9,-2 \left(\sqrt{2}+2 \sqrt{4-2 \sqrt{2}}-4\right),\\
&\quad 8 \sqrt{2}-6 \sqrt{2 \left(\sqrt{2}+2\right)}+5,2 \left(-2 \sqrt{2}+\sqrt{14 \sqrt{2}+20}-3\right)\} \\
&{\approx \{1., 0., 1.05595, 0.842004, 0.634953, 0.960434\}},\\
{\vec{\xi}}_2 &= \{1,0,2 \sqrt{2-\sqrt{2}}-\frac{3}{\sqrt{2}}+2,2 \sqrt{\sqrt{2}+2}-\frac{1}{\sqrt{2}}-2, \\
&\quad -2 \left(\sqrt{2}+2 \sqrt{4-2 \sqrt{2}}-4\right),6 \sqrt{2}+4 \sqrt{10-7 \sqrt{2}}-9\} \\
&{\approx \{1., 0., 1.40941, 0.988411, 0.842004, 0.753383\},}\\
{\vec{\xi}}_3 & \approx \{1,0.447769,0,1.21903,-0.161745,2.4715\},
\end{aligned} \]
and ${\vec{\xi}}_{3+i}$ is obtained from ${\vec{\xi}}_i$ via \eqref{perm}, $i=1,2,3$.
 
If this is the case, then in fact $C(A)$ is given by \eqref{ee0}, with $E_0$ being an ellipse with the foci $\pm\sqrt{2-\sqrt{2}}$. Moreover, the half minor axes of $\pm E$ and $E_0$ are $\left\{1,\sqrt{9 \sqrt{2}-4 \sqrt{14 \sqrt{2}+20}+15}\right\}t$, $\left\{1,\sqrt{2 \left(\sqrt{2}+\sqrt{4-2 \sqrt{2}}-1\right)}\right\}t$, $\left\{1.20323,1.44257\right\}t$, respectively. 
\end{thm} 

The proof of Theorem~\ref{th:7outer} follows the same logic as that of Theorem~\ref{th:7inner}.  Note that in the subcase $\xi_3 = 0$ there is only one one-parametric family of solutions. This explains the difference between the number of ${\vec{\xi}}_i$ in Theorems~\ref{th:7inner} and \ref{th:7outer}. 

The Kippenhahn curve of the matrix $A$ corresponding to $\xi_j$ given by $\vec{\xi}_1$ and $t = 1$ is shown below.
\begin{figure}[H]
    \centering
    \includegraphics[width=0.5\linewidth]{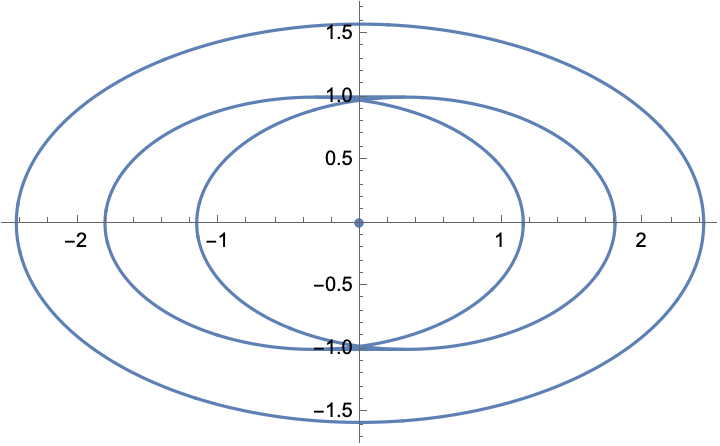}
    \begin{minipage}{1.0\textwidth}
        \centering
        \caption{
            $\vec{\xi} \to \{1,0,7 \sqrt{2}-4 \sqrt{7 \sqrt{2}+10}+9,-2 \left(\sqrt{2}+2 \sqrt{4-2 \sqrt{2}}-4\right), \quad 8 \sqrt{2}-6 \sqrt{2 \left(\sqrt{2}+2\right)}+5,2 \left(-2 \sqrt{2}+\sqrt{14 \sqrt{2}+20}-3\right)\}$}
    \end{minipage}
\end{figure}
\providecommand{\bysame}{\leavevmode\hbox to3em{\hrulefill}\thinspace}
\providecommand{\MR}{\relax\ifhmode\unskip\space\fi MR }
\providecommand{\MRhref}[2]{%
  \href{http://www.ams.org/mathscinet-getitem?mr=#1}{#2}
}
\providecommand{\href}[2]{#2}

\end{document}